\documentclass[12pt]{article}

\usepackage[german, american]{babel}  

\usepackage{times}
\usepackage{graphicx}
\usepackage{amsfonts}       %
\usepackage{amsmath} 
\usepackage{amssymb}        %
\usepackage{amstext}        %
\usepackage{amsthm}         

\usepackage{epsfig}
\psfigdriver{dvips}
\usepackage{latexsym}
\usepackage[matrix,curve]{xypic}
\usepackage{rotating}
\rotdriver{dvips}

\newcommand{\Z}{\mathbb{Z}}
\newcommand{\Q}{\mathbb{Q}}

\DeclareMathSizes{12}{11}{8}{5} 
\usepackage{xypic,dsfont}

\begin{document}
\author{Rainer Weissauer}
\title{Tannakian Categories attached to abelian Varieties}
\maketitle


\bigskip\noindent

\bigskip
Let $k$ be an algebraically closed field $k$, where $k$ is either the algebraic closure of a
finite field or a field of characteristic zero. Let $l$ be a prime different from the
characteristic of $k$.

\bigskip\noindent
{\it Notations.} For a variety $X$ over $k$ let $D_c^b(X,\overline\Q_l)$ denote the
triangulated category of complexes of etale $\overline \Q_l$-sheaves on $X$ in the sense of
\cite{KW}. For a complex $K\in D_c^b(X,\overline\Q_l)$ let $D(K)$ denote its Verdier dual, and
${\cal H}^\nu(K)$ denote its etale cohomology $\overline \Q_l$-sheaves with respect to the
standard $t$-structure. The abelian subcategory $Perv(X)$ of middle perverse sheaves is the
full subcategory of all $K\in D_c^b(X,\overline\Q_l)$, for which  $K$ and its Verdier dual
$D(K)$ are contained in the full subcategory ${}^p\!D^{\leq 0}(X)$ of semi-perverse sheaves,
where $L\in D_c^b(X,\overline\Q_l)$ is semi-perverse if and only if $dim(S_\nu) \leq \nu$ holds
for all integers $\nu \in \Z$, where $S_\nu$ denotes the support \label{supports} of the
cohomology sheaf ${\cal H}^{-\nu}(L)$ of $L$.

\bigskip
If $k$ is the algebraic closure of a finite field $\kappa$, then a complex $K$ of etale
$\overline\Q_l$-Weil sheaves is mixed of weight $\leq w$, if all its cohomology sheaves ${\cal
H}^\nu(K)$ are mixed etale $\overline\Q_l$-sheaves with upper weights $w({\cal H}^\nu(K))-\nu
\leq w$ for all integers $\nu$. It is called pure of weight $w$, if $K$ and its Verdier dual
$D(K)$ are mixed of weight $\leq w$.  Concerning base fields of characteristic zero, we assume
mixed sheaves to be sheaves of geometric origin in the sense of the last chapter of \cite{BBD},
so we still dispose over the notion of the weight filtration and purity and Gabber's
decomposition theorem in this case. In this sense let $Perv_m(X)$ denote the abelian category
of mixed perverse sheaves on $X$. The full subcategory $P(X)$ of $Perv_m(X)$ of pure perverse
sheaves is a semisimple abelian category.

\bigskip\noindent
{\it Abelian varieties.} Let $X$ be an abelian variety $X$ of dimension $g$ over an
algebraically closed field $k$. The addition law of the abelian variety $a:X\times X\to X$
defines the convolution product $K*L\in D_c^b(X,\overline\Q_l)$ of two complexes $K$ and $L$ in
$D_c^b(X,\overline\Q_l)$ by the direct image
$$  K*L \ = \ Ra_*(K\boxtimes L) \ .$$

\bigskip\noindent
For the skyscraper sheaf $\delta_0$ concentrated at the zero element $0$ notice $ K*\delta_0 =
K$.

\goodbreak
\bigskip\noindent
{\it Translation-invariant sheaf complexes}. More generally $K*\delta_{x} = T^*_{-x}(K)$, where
$x$ is a closed $k$-valued point in $X$, $\delta_x$ the skyscraper sheaf with support in
$\{x\}$ and where $T_{x}(y)=y+x$ denotes the translation $T_x:X\to X$ by $x$. In fact $
T^*_y(K*L) \cong T^*_y(K)*L \cong K*T^*_y(L) $ holds for all $y\in X(k) $.
 For $K\in D_c^b(X,\overline\Q_l)$ let $Aut(K)$ be the abstract group of all closed $k$-valued points $x$ of
$X$, for which $T^*_x(K)\cong K$ holds.
 A complex $K$ is called translation-invariant, provided $Aut(K)=X(k)$.
If $f:X\to Y$ is a surjective homomorphism between abelian varieties, then the direct image
$Rf_*(K)$ of a translation-invariant complex is translation-invariant.  As a consequence of the
formulas above, the convolution of an arbitrary $K\in D_c^b(X,\overline\Q_l)$ with a
translation-invariant complex on $X$ is a translation-invariant complex.  A
translation-invariant perverse sheaf $K$ on $X$ is of the form $K=E[g]$, for an ordinary etale
translation-invariant $\overline\Q_l$-sheaf $E$. For a translation-invariant complex $K\in
D_c^b(X,\overline\Q_l)$ the irreducible constituents of the perverse cohomology sheaves ${}^p
H^\nu(K)$ are translation-invariant.

\bigskip\noindent
{\it Multipliers}. The subcategory $T(X)$ of $Perv(X)$ of all perverse sheaves, whose
irreducible perverse constituents are translation-invariant, is a Serre subcategory of the
abelian category $Perv(X)$.  Let denote $ \overline{Perv}(X)$ its abelian quotient category and
$\overline P(X)$ the image of $P(X)$, which is a full subcategory of semisimple objects. The
full subcategory of $D_c^b(X,\overline \Q_l)$ of all $K$, for which ${}^pH^\nu(K)\in T(X)$, is
a thick subcategory of the triangulated category $D_c^b(X,\overline\Q_l)$. Let
$$ \overline{D}_c^b(X,\overline\Q_l)$$
be the corresponding triangulated quotient category, which contains $ \overline{Perv}(X)$. Then
the convolution product
$$ *:\ \overline D_c^b(X,\overline\Q_l) \times \overline D_c^b(X,\overline\Q_l) \to \overline
D_c^b(X,\overline\Q_l) $$ still is well defined, by reasons indicated above. 

\bigskip\noindent
{\bf Definition.} {\it A perverse sheaf $K$ on $X$ is called a multiplier, if the convolution
induced by  $K$
$$ *K: \overline {D_c^b}(X,\overline{\Q_l}) \to \overline {D_c^b}(X,\overline{\Q_l}) $$
preserves the abelian subcategory $\overline{Perv}(X)$.}

\bigskip
Obvious from this definition are the following properties of multipliers: If $K$ and $L$ are
multipliers, so are the product $K*L$ and the direct sum $K\oplus L$. Direct summands of
multipliers are multipliers. If $K$ is a multiplier, then the Verdier dual $D(K)$ is a
multiplier and also the dual
$$ K^\vee = (-id_X)^*(D(K)) \ .$$

\bigskip\noindent
\underbar{Examples}: 1) Skyscraper sheaves are multipliers 2) If $i:C\hookrightarrow X$ is a
projective curve, which generates the abelian variety $X$, and $E$ is an etale
$\overline\Q_l$-sheaf on $C$ with finite monodromy, then the intersection cohomology sheaf
attached to $(C,E)$ is a multiplier. 3) If $:Y \hookrightarrow X$ is a smooth ample divisor,
then the intersection cohomology sheaf of $Y$ is a multiplier.

\bigskip\noindent
{\it The proofs}. 1) is obvious. For 2) we gave in \cite{BN} a proof by reduction mod $p$ using
the Cebotarev density theorem and counting of points. Concerning 3) the morphism $j:
U=X\setminus Y \hookrightarrow X$ is affine for ample divisors $Y$. Hence
$\lambda_U=Rj_!\overline\Q_l[g]$ and $\lambda_Y=i_* \overline\Q_{l,Y}[g-1]$ are perverse
sheaves, which coincide in $\overline{Perv}(X)$. The morphism $\pi = a\circ (j\times id_X)$ is
affine. Indeed $W=\pi^{-1}(V)$ is affine for affine subsets $V$ of $X$, $W$ being isomorphic
under the isomorphism $(u,v)\mapsto (u,u+v)$ of $X^2$ to the affine product $U\times V$. By the
affine vanishing theorem of Artin: For perverse sheaves $L\in Perv(X)$ we get $\lambda_U
\boxtimes L \in Perv(X^2)$ and $ {}^p H^\nu(R\pi_! (\lambda_U \boxtimes L)) = 0$ for all $\nu <
0 $. The distinguished triangle $ \bigl(Ra_*(\lambda_Y\boxtimes L) , R\pi_!(\lambda_U \boxtimes
L) , Ra_*(\delta_X \boxtimes L)\bigr) $ for $\delta_X=\overline\Q_{l,X}[g]$ and the
corresponding long exact perverse cohomology sequence gives isomorphisms $ {}^p H^{\nu
-1}(\delta_X * L) \cong {}^p H^{\nu}(\lambda_Y * L)$ for the integers $\nu <0$. Since
$Ra_*(\delta_X \boxtimes L) = \delta_X * L$ is a direct sum of translates of constant perverse
sheaves $\delta_X$, we conclude $ {}^p H^{\nu}(\lambda_Y
* L)$ for $\nu<0 $ to be zero in $\overline{Perv}(X)$.
For smooth $Y$  the intersection cohomology sheaf is $\lambda_Y=i_* \overline\Q_{l,Y}[g-1]$,
and it is self dual. Hence by Verdier duality $i_* \overline\Q_{l,Y}[g-1] * L$  has image in
$\overline{Perv}(X)$. Thus $i_* \overline\Q_{l,Y}[g-1]$ is a multiplier.  $\quad  \square$


\goodbreak
\bigskip
Let $M(X)\subseteq P(X)$ denote the full category of semisimple multipliers. Let ${\overline
M}(X)$ denote its image in the quotient category $\overline{P}(X)$ of $P(X)$. Then, by the
definition of multipliers, the convolution product preserves ${\overline M}(X)$
$$ *:\ {\overline M}(X) \times {\overline M}(X) \to {\overline M}(X)
\ .$$

\bigskip\noindent
{\bf Theorem}. {\it With respect to this convolution product the category ${\overline M}(X)$ is
a semisimple super-Tannakian $\overline \Q_l$-linear tensor category, hence as a tensor
category ${\overline M}(X)$ is equivalent to the category of representations
$Rep(G,\varepsilon)$ of a projective limit $$G = G(X) $$ of supergroups.}

\goodbreak
\bigskip\noindent
{\it Outline of proof}. The convolution product obviously satisfies the usual commutativity and
associativity constraints compatible with unit objects. See \cite{BN} 2.1. By \cite{BN},
corollary 3 furthermore one has functorial isomorphisms
$$ Hom_{{\overline M}(X)}(K,L) \ \cong \ \Gamma_{\{0\}}(X,{\cal
H}^0(K*L^\vee)^*) \ ,$$ where ${\cal H}^0$ denotes the degree zero cohomology sheaf and
$\Gamma_{\{0\}}(X,-)$ sections with support in the neutral element. Let $L=K$ be simple and
nonzero. Then the left side becomes $End_{{\overline M}(X)}(K) \cong \overline{\Q_l}$. On the
other hand $K*L^\vee$ is a direct sum of a perverse sheaf $P$ and translates of
translation-invariant perverse sheaves. Hence ${\cal H}^0(K*L^\vee)^\vee)$ is the direct sum of
a skyscraper sheaf $S$ and translation-invariant etale sheaves. Therefore
$\Gamma_{\{0\}}(X,{\cal H}^0(K*L^\vee)^\vee) = \Gamma_{\{0\}}(X,S)$. By a comparison of both
sides therefore $S=\delta_{0}$. Notice $\delta_0$ is the unit element $1$ of the convolution
product. Using the formula above we not only get $$ Hom_{{\overline M}(X)}(K,L) \ \cong \
Hom_{{\overline M}(X)}(K*L^\vee,1) \ ,
$$ but also find a nontrivial morphism $$ ev_K: K*K^\vee \to 1 \ .$$ By semisimplicity
$\delta_{0}$ is a direct summand of the complex $K*K^\vee$. In particular the K\"{u}nneth
formula implies, that the etale cohomology groups  do not all vanish identically
$$ H^\bullet(X,K)\neq 0 \ .$$ Therefore the arguments of \cite{BN} 2.6 show, that the simple perverse
sheaf $K$ is dualizable. Hence ${\overline M}(X)$ is a rigid $\overline{\Q_l}$-linear tensor
category. Let ${\cal T}$ be a finitely $\otimes$-generated tensor subcategory with generator
say $A$. To show ${\cal T}$ is super-Tannakian, by \cite{De} it is enough to show for all $n$
$$ lenght_{\cal T}(A^{*n}) \leq N^n \ ,$$
where $N$ is a suitable constant. For any $B\in {\overline M}(X)$ let $B$, by abuse of
notation, also denote the perverse semisimple representative in $Perv(X)$ without translation
invariant summand. Put $h(B,t)= \sum_\nu dim_{\overline{\Q_l}}(H^\nu(X,B))t^\nu$. Then $
lenght_{\cal T}(B) \leq h(B,1)$, since every summand of $B$ is a multiplier and therefore has
nonvanishing cohomology. For $B=A^{*n}$ the K\"{u}nneth formula gives $h(B,1) = h(A,1)^n$.
Therefore the estimate above holds for $N=h(A,1)$. This completes the outline for the proof of
the theorem. $\quad \square$

\bigskip\noindent
{\it Principally polarized abelian varieties}. Suppose $Y$ is a divisor in $X$ defining a
principal polarization. Suppose the intersection cohomology sheaf $\delta_Y$ of $Y$ is a
multiplier. Then a suitable translate of $Y$ is symmetric, and again a multiplier. So we may
assume $Y=-Y$ is symmetric. Let ${\overline M}(X,Y)$ denote the super-Tannakian subcategory of
$\overline M(X)$ generated by $\delta_Y$. The corresponding super-group $G(X,Y)$ attached to
${\overline M}(X,Y)$ acts on the super-space $W=\omega(\delta_Y)$ defined by the underlying
super-fiber functor $\omega$ of $\overline M(X)$. By assumption $\delta_Y$ is self dual in the
sense, that there exists an isomorphism $\varphi: \delta_Y^\vee \cong \delta_Y$. Obviously
$\varphi^\vee = \pm \varphi$. This defines a nondegenerate pairing on $W$, and the action of
$G(X,Y)$ on $W$ respects this pairing.

\bigskip\noindent
{\it Curves}. If $X$ is the Jacobian of smooth projective curve $C$ of genus $g$ over $k$, $X$
carries a natural principal polarization $Y=W_{g-1}$. If we replace this divisor by a symmetric
translate, then $Y$ is a multiplier. The corresponding group $G(X,Y)$ is the semisimple
algebraic group $G=Sp(2g-2,\overline\Q_l)/\mu_{g-1}[2]$ or $G=Sl(2g-2,\overline\Q_l)/\mu_{g-1}$
depending on whether the curve $C$ is hyperelliptic or not. The representation $W$ of $G(X,Y)$
defined by $\delta_Y$ as above is the unique irreducible $\overline\Q_l$-representation of
$G(X,Y)$ of highest weight, which occurs in the $(g-1)$-th exterior power of the
$(2g-2)$-dimensional standard representation of $G$. See \cite{BN}, section 7.6.

\bigskip\noindent
{\it Conjecture}. One could expect, that a principal polarized abelian variety $(X,Y)$ of
dimension $g$ is isomorphic to a Jacobian variety $(Jac(C),W_{g-1})$ of a smooth projective
curve $C$ (up to translates of the divisor $Y$ in $X$ as explained above) if and only if $Y$ is
a multiplier with corresponding super-Tannakian group $G(X,Y)$ equal to one of the two groups
$$Sp(2g-2,\overline\Q_l)/\mu_{g-1}[2] \ \mbox{  or }\ Sl(2g-2,\overline\Q_l)/\mu_{g-1}\ .$$

\end{document}